\documentclass[12pt]{amsart}
\usepackage{amsmath,amssymb,amsbsy,amsfonts,amsthm,latexsym,
                       amsopn,amstext,amsxtra,euscript,amscd}
\usepackage{enumerate}
\usepackage{array}
\usepackage{ifthen}
\usepackage{url}
\usepackage{color}

\newtheorem{theorem}{Theorem}
\newtheorem{lemma}{Lemma}

\newtheorem*{corollary}{Corollary}

\theoremstyle{definition}

\newcommand{\klammern}[4][]%
{\ifthenelse{\equal{#1}{}}{\left#2}{\csname#1\endcsname#2}%
#4\ifthenelse{\equal{#1}{}}{\right#3}{\csname#1\endcsname#3}}

\def\KK{\mathbb K}
\def\QQ{\mathbb Q}
\def\ZZ{\mathbb Z}

\begin{document}

\title[$P_n$ or $Q_n$ which are concatenations of two repdigits in base $b$ ]{\small Pell or Pell-Lucas numbers as concatenations of two repdigits in base $b$ }

\author[K. N. Ad\'edji, A. Filipin, S. E. Rihane and A. Togb\'e]{Kou\`essi Norbert Ad\'edji, Alan Filipin, Salah Eddine Rihane and Alain Togb\'e}

\keywords{Pell numbers, Pell-Lucas numbers, $b$-repdigits, linear forms in logarithms, Diophantine equations, Reduction method}

\begin{abstract}
Let $b$ be a positive integer such that $2\le b\le 10.$ In this study, we find all Pell or Pell-Lucas numbers as concatenations of two repdigits in base $b.$ As a corollary, it is show that the largest Pell or Pell-Lucas numbers which can be expressible as a concatenations of two repdigits in base $b$ are $P_{11} = 5741$ and $Q_{5}=82,$ respectively.
\end{abstract}

\date{\today}

\maketitle

\section{Introduction}\label{sec:1}

In the literature, especially in mathematics and physics, there are a lot of integer sequences, which are used in almost every field modern sciences. Admittedly, the Fibonacci sequence is one of the most famous and curious numerical sequences in mathematics and have been widely studied from both algebraic and combinatorial prospectives. Also, there is the Pell sequence, which is as important as the Fibonacci sequence. The Pell sequence $\{P_n\}$ are defined by recurrence $P_n =2P_{n-1}+ P_{n-2},\; n\ge 2$ with $P_0=0$ and $P_1=1$ and the Pell-Lucas sequence $\{Q_n\}$  by the same recurrence but with initial conditions $Q_0=Q_1=2.$ The explicit Binet formulas for $\{P_n\}$ and $\{Q_n\}$ are
\begin{align}
P_n=\dfrac{\alpha^n-\beta^n}{\alpha-\beta}\quad \text{and}\quad Q_n=\alpha^n+\beta^n\quad \text{for~all}\; n\ge 0,
\end{align}
where $\alpha=1+\sqrt{2}$ and $\beta=1-\sqrt{2}$ are the roots of the characteristic equation $x^2-2x-1 =0.$ It can be seen that $2 <\alpha< 3$ and $-1<\beta<0.$ The inequalities 
\begin{align}\label{P_n-Q_n}
\alpha^{n-2}\le P_n\le \alpha^{n-1}\quad\text{and}\quad \alpha^{n-2}\le Q_n< \alpha^{n+1},
\end{align}
are well known, where $n\ge 1.$ Given an integer $b>1,$ a base $b$-repdigit is a number $N$ of the form
\[
N=\dfrac{d(b^m-1)}{b-1}=\overline{\underbrace{d\ldots d}_{m~times}}~_{(b),}
\]
for some positive integers $d,m$ with $d\in \{1,\ldots,b-1\}$ and $m\ge 1.$ In this work, we study the problem of finding all Pell and Pell-Lucas numbers that are concatenations of two $b$-repdigits where $2\le b\le 10$. More precisely, we completely solve the following two Diophantine equations
\begin{align}\label{main-equation_1}
P_n=\overline{\underbrace{d_1\ldots d_1}_
{\ell_1~times}\underbrace{d_2\ldots d_2}_
{\ell_2~times}}~_{(b)}=d_1\left(\dfrac{b^{\ell_1}-1}{b-1}\right)\times b^{\ell_2}+d_2\left(\dfrac{b^{\ell_2}-1}{b-1} \right),
\end{align}
and
\begin{align}\label{main-equation_2}
 Q_n=\overline{\underbrace{d_1\ldots d_1}_
{\ell_1~times}\underbrace{d_2\ldots d_2}_
{\ell_2~times}}~_{(b)}=d_1\left(\dfrac{b^{\ell_1}-1}{b-1}\right)\times b^{\ell_2}+d_2\left(\dfrac{b^{\ell_2}-1}{b-1} \right),
\end{align}
in non-negative integers $(n, d_1, d_2, \ell_1,\ell_2)$ with $n\ge 1,\; d_1\ne d_2$, and $d_1,d_2\in \{1,\ldots,b-1\}.$

In view of the above, our main results of this paper is as follows.
\begin{theorem}\label{main:1}
The only Pell numbers that are concatenations of two repdigits in base $b$ with $2\le b\le 10$ are
\[
2, 5, 12, 29, 70, 169, 408, 5741.
\]
More precisely, we have
$$
\begin{array}{ccccl}
2&=&P_2&=&\overline{10}_2,\\
5&=&P_3&=&\overline{12}_3=\overline{10}_5,\\
12&=&P_4&=&\overline{1100}_2=\overline{110}_3=\overline{30}_4=\overline{15}_7=\overline{14}_8=\overline{13}_9=\overline{12}_{10},\\
29&=&P_5&=&\overline{45}_6=\overline{41}_7=\overline{35}_8=\overline{32}_9=\overline{29}_{10},\\
70&=&P_6&=&\overline{70}_{10},\\
169&=&P_7&=&\overline{2221}_4=\overline{441}_6=\overline{331}_7,\\
408&=&P_8&=&\overline{1122}_7,\\
5741&=&P_{11}&=&\overline{7778}_9.
\end{array}
$$
\end{theorem}
\begin{theorem}\label{main:2}
The only Pell-Lucas numbers that are concatenations of two repdigits in base $b$ with $2\le b\le 10$ are
\[
2, 6, 14, 34, 82.
\]
More precisely, we have
$$
\begin{array}{ccccl}
2&=&Q_0&=&\overline{10}_2,\\
2&=&Q_1&=&\overline{10}_2,\\
6&=&Q_2&=&\overline{110}_2=\overline{20}_3=\overline{12}_4=\overline{10}_6,\\
14&=&Q_3&=&\overline{1110}_2=\overline{112}_3=\overline{32}_4=\overline{24}_5=\overline{20}_7=\overline{16}_8=\overline{15}_9=\overline{14}_{10},\\
34&=&Q_4&=&\overline{114}_5=\overline{54}_6=\overline{46}_7=\overline{42}_8=\overline{37}_9=\overline{34}_{10}\\
82&=&Q_5&=&\overline{122}_8=\overline{82}_{10}.
\end{array}
$$
\end{theorem}
The two theorems above allow us to deduce the following result
\begin{corollary}
The largest Pell and Pell-Lucas numbers which can be representable as a concatenations of two repdigits in base $b$ are 
$P_{11}=5741=\overline{7778}_{9}$ and $Q_{5} = 82=\overline{122}_{8}=\overline{82}_{10},$ respectively.
    \end{corollary}
This paper is inspired by the result of Alahmadi, Altassan, Luca, and Shoaib \cite{AALS:2019}, in which they find all Fibonacci numbers that are concatenations of two repdigits. Our method of proof involves the application of Baker's theory for linear forms in logarithms of algebraic numbers, and the Baker-Davenport reduction procedure. Computations are done with the help of a computer program in Maple. The outline for this article is as follows. In Section~\ref{sec:2} we will list the main results that we will use in order to establish Theorems~\ref{main:1} and \ref{main:2}. Finally, in Sections~\ref{sec:3} and \ref{sec:4} we will prove Theorems~\ref{main:1} and \ref{main:2} respectively.

\section{Some useful results}\label{sec:2}

To solve the Diophantine equations involving repdigits and the terms of binary recurrence sequences, many authors have used Baker's theory to reduce lower bounds concerning linear forms in logarithms of algebraic numbers. These lower bounds play an important role while solving such Diophantine  equation. We start with recalling some basic definitions and results from algebraic number theory. For any non-zero algebraic $\eta$ of degree $d$ over $\QQ$, whose minimal polynomial over $\ZZ$ is $a_0 \prod_{j=1}^{d}(X-\eta^{(j)}),$ we denote by 
$$
 h(\eta)= \frac{1}{d}\left(\log |a_0|+\sum_{j=1}^{d}\log\max\left(1,\left|\alpha^{(j)} \right| \right) \right)
$$ 
its absolute logarithmic height. Note that, if $\eta=\frac{p}{q}\in \QQ$ is a rational number in reduced form with $q>0$, then the above definition reduces to $h(\eta)=\log\max\{|p|, q\}.$ We list some well known properties of the height function below, which we shall subsequently use without reference:
\begin{eqnarray}
\label{equt_10}h(\eta_1 \pm \eta_2) &\leq& h(\eta_1)+ h(\eta_2) +\log2,\\
 \label{equt_11}h(\eta_1\eta_2^{\pm}) &\leq& h(\eta_1) + h(\eta_2),\\
\label{equt_12}h(\eta^s)&=&|s|h(\eta),\quad (s\in \ZZ).
\end{eqnarray}
We quote the version of Baker's theorem proved by Bugeaud, Mignotte and Siksek in (\cite{BMS:2006}, Theorem~9.4, pp. 989).

\begin{lemma}[Theorem~2 of \cite{BGL:2016}]\label{Matveev}
Assume that $ \gamma_1,\ldots,\gamma_t$ are positive real algebraic numbers in a real algebraic number field $\KK$ of degree $D,$ $b_1,\ldots,b_t$ are rational integers, and 
$$  
\Lambda:=\gamma_1^{b_1}\cdots\gamma_t^{b_t}-1,
$$
is not zero. Then 
$$  
|\Lambda|>\exp\left(-1.4\cdot30^{t+3}\cdot t^{4.5}\cdot D^2(1+\log D)(1+\log B)A_1\cdots A_t\right),
$$
where 
$$  
B\ge \max\{|b_1|,\ldots,|b_t|\},
$$
and 
$$
A_i\ge \max\{Dh(\gamma_i),|\log\gamma_i|,0.16\},\quad \mbox{for all}\quad i=1,\ldots,t.
$$
\end{lemma}
Using the above theorem and properties of logarithmic height, we will obtain upper bounds for the $b$-repdigits and the index of Pell or Pell-Lucas numbers. Then, we will apply the following lemma for the further reduction of the obtained upper bounds so that in the remaining range the Pell and Pell-Lucas numbers which are concatenations of two repdigits in base $b$ can be verified with direct computation.
\begin{lemma}[Lemma~5a of \cite{Dujella-Peto}]\label{Dujella-Petho}
Let $M$ be a positive integer, $p/q$ be a convergent of the continued fraction expansion of the irrational number $\tau$ such that $q>6M,$ and $A, B, \mu$ be some real numbers with $A>0$ and $B>1.$ Furthermore, let
$$ 
\varepsilon:=\parallel\mu q\parallel-M\cdot\parallel\tau q\parallel.
$$
 If $\varepsilon>0,$ then there is no solution to the inequality
\begin{equation}
0<|u\tau-v+\mu|<AB^{-w}
\end{equation}
in positive integers $u,v$ and $w$ with
$$
 u\leq M\; \mbox{and}\; w\geq \frac{\log(Aq/\varepsilon)}{\log B}.
$$
\end{lemma}
We will use the following well known properties of continud fraction.
\begin{lemma}[Pages $30$ and $37$ in \cite{Khinchin:1963}]\label{Legendre}
Let $\tau$ be an irrational number.
\begin{itemize}
\item[(i)] If $x,y$ are positive integers such that
$$
\left|\tau-\frac{x}{y}\right|<\frac{1}{2y^2},
$$
then $x/y=p_k/q_k$ is a convergent of $\tau$.
\item[(ii)] Let $M$ be a positive real number and $p_0/q_0,p_1/q_1,\ldots$ be all the convergents of the continued fraction of $\tau$. Let $N$ be the smallest positive  integer such that $q_N>M$. Put
$a(M):=\max\{a_k: k=0,1,\ldots,N\}.$ Then, the inequality
$$
\left|\tau-\dfrac{x}{y}\right|>\dfrac{1}{(a(M)+2)y^2},
$$
holds for all pairs $(x,y)$ of integers with $0<y<M$.
\end{itemize}
\end{lemma}

\section{Proof of Theorem~\ref{main:1}}\label{sec:3}

\subsection{Bounding $n$ and $l_1+l_2$}

To start with, consider the following Diophantine equation which is equivalent to \eqref{main-equation_1}
\begin{align}\label{P_n=rep}
P_n=\dfrac{1}{b-1}\left( d_1 b^{l_1+l_2}-(d_1-d_2)b^{l_2}-d_2\right),
\end{align}
where $d_1, d_2, l_1$ and $l_2$ are non negative integers with $d_1, d_2\in \{1,\ldots,b-1\}$ and $2\le b\le 10.$ We assume that $n>110.$ From inequality \eqref{P_n-Q_n}, we can get 
$$
\alpha^{n-2}\le P_n<b^{l_1+l_2}\quad \text{and}\quad b^{l_1+l_2-1}<P_n\le \alpha^{n-1},
$$
which implies that
\begin{align*}
(l_1+l_2)\log b+\log\alpha -\log b<n\log\alpha<(l_1+l_2)\log b+2\log \alpha.
\end{align*}
Since $\log\alpha-\log 10<\log\alpha-\log b,$ we can get 

\begin{align}\label{nlog-alpha}
(l_1+l_2)\log b-1.5<n\log\alpha<(l_1+l_2)\log b+1.8.
\end{align}
From \eqref{nlog-alpha}, we get 
$$
l_1+l_2>\dfrac{n\log\alpha-1.8}{\log 10}>41.
$$
Using  \eqref{P_n=rep} and Binet's formula for Pell sequence, we get
\[
(b-1)\alpha^n-2\sqrt{2}d_1b^{l_1+l_2}=(b-1)\beta^n-2\sqrt{2}\left[(d_1-d_2)b^{l_2}+d_2\right].
\]
Notice that $|d_1-d_2|\le b-2<8.$ Since $n>110$, we have 
\begin{align*}
\left| (b-1)\alpha^n-2\sqrt{2}d_1b^{l_1+l_2}\right|&=\left|(b-1)\beta^n-2\sqrt{2}\left[(d_1-d_2)b^{l_2}+d_2\right] \right|\\
&\le 9\alpha^{-n}+2\sqrt{2}(8b^{l_2}+9)\\
&<25.5\cdot b^{l_2},
\end{align*}
which implies that
\begin{align}\label{Gamma_1-bound}
\left|\dfrac{b-1}{2\sqrt{2}d_1}\cdot \alpha^n\cdot b^{-(l_1+l_2)} -1\right|<\dfrac{9.1}{b^{l_1}}.
\end{align}
Let 
\begin{align}\label{Gamm_1}
\Gamma_1:=\dfrac{b-1}{2\sqrt{2}d_1}\cdot \alpha^n\cdot b^{-(l_1+l_2)} -1.
\end{align}
It is easy to see that $\Gamma_1\neq 0.$ Indeed if $\Gamma_1=0,$ then 
$$
\alpha^{2n}=\dfrac{8d_1^2b^{2(l_1+l_2)}}{(b-1)^2},
$$
which is a contradiction since $\alpha^{2n}$ is irrational for $n\ge 1.$  According to Lemma~\ref{Matveev} we can take $t=3$ and 
\[
(\gamma_1, b_1):=\left(\dfrac{b-1}{2\sqrt{2}d_1}, 1 \right),\quad (\gamma_2, b_2):=(\alpha, n),\quad (\gamma_3, b_3):=(b, -l_1-l_2).
\]
Thus, we have $\KK=\QQ(\gamma_1,\gamma_2,\gamma_3)=\QQ(\alpha)$, $D=[\KK:\QQ]=2.$ Based on the inequality 
$$
b^{l_1+l_2-1}<P_n\le \alpha^{n-1},
$$ 
we deduce that
$$
l_1+l_2<n\dfrac{\log\alpha}{\log 2}+\dfrac{\log(2/\alpha)}{\log 2}<1.3n.
$$
As $B\ge \max\{|1|, |n|, |-l_1-l_2|\},$ we can take $B:=1.3n.$
Note that
\begin{align*}
h(\gamma_1)&=h\left(\dfrac{b-1}{2\sqrt{2}d_1} \right)\\
&\le h\left(\dfrac{b-1}{2d_1}\right)+h(\sqrt{2})=\log (\max\{b-1, 2d_1\})+\dfrac{1}{2}\log 2\\
&\le 2\log 9+\dfrac{1}{2}\log 2.
\end{align*}
Moreover, $h(\gamma_2)=h(\alpha)=\dfrac{1}{2}\log\alpha$ and $h(\gamma_3)=\log b\le \log 10.$ Thus, we can take
\[
A_1:=9.5,\quad A_2:= 0.89\quad \text{and}\quad A_3:=4.7.
\]
Hence, the Lemma~\ref{Matveev} allows us to obtain
\begin{align}\label{log-Gamma_1}
\log|\Gamma_1|>-C_1(1+\log 1.3n),
\end{align}
where $C_1=3.853\times 10^{13}.$
Thus from \eqref{Gamma_1-bound} and \eqref{log-Gamma_1}, we can get
\begin{align}\label{l_1log b}
l_1\log b<C_1(1+\log 1.3n)+\log 9.1.
\end{align}
We rewrite equation \eqref{P_n=rep}, then we get
\[
\alpha^n-2\sqrt{2}\left(\dfrac{d_1b^{l_1}-(d_1-d_2)}{b-1} \right)b^{l_2}=\beta^n-\dfrac{2\sqrt{2} d_2}{b-1}.
\]
Since $n>110,$ we deduce that
\[
\left|\alpha^n-2\sqrt{2}\left(\dfrac{d_1b^{l_1}-(d_1-d_2)}{b-1} \right)b^{l_2}  \right|=\left| \beta^n-\dfrac{2d_2\sqrt{2}}{b-1} \right|\le \alpha^{-n}+2\sqrt{2}<3.
\]
It follows that
\begin{align}\label{Gamma_2-bound}
\left|\left(\dfrac{2(d_1b^{l_1}-(d_1-d_2))\sqrt{2}}{b-1} \right)\cdot\alpha^{-n}\cdot b^{l_2}-1\right|<\dfrac{3}{\alpha^n}.
\end{align}
Let 
\begin{align}\label{Gamma_2}
\Gamma_2:=\left(\dfrac{2(d_1b^{l_1}-(d_1-d_2))\sqrt{2}}{b-1} \right)\cdot\alpha^{-n}\cdot b^{l_2}-1,
\end{align}
then $\Gamma_2\neq 0.$ If we assume that $\Gamma_2=0,$ then we get the following equation
\[
\alpha^{2n}=\dfrac{8(d_1b^{l_1}-(d_1-d_2))^2}{(b-1)^2}b^{2l_2}\in \QQ,
\]
which is impossible for $n\ge 1.$ According to Lemma~\ref{Matveev} and using \eqref{Gamma_2}, we can take the following data
\begin{align*}
t:=3,\quad \gamma_1:=\dfrac{2(d_1b^{l_1}-(d_1-d_2))\sqrt{2}}{b-1},\quad \gamma_2:=\alpha,\quad \gamma_3:=b,
\end{align*}
and the exponents 
$$
b_1:=1,\quad b_2:=-n,\quad b_3:=l_2.
$$
Thus, we have $\KK=\QQ(\gamma_1,\gamma_2,\gamma_3)=\QQ(\alpha)$, $D=[\KK:\QQ]=2.$ Since $B\ge \{|1|, |-n|,|l_2|\}$ and we also knew that $l_1+l_2<1.3n,$ so we can take $B=1.3n.$ From \eqref{l_1log b}, we can get
\[
\begin{array}{lcl}
h(\gamma_1)&=& h\left( \dfrac{2(d_1b^{l_1}-(d_1-d_2))\sqrt{2}}{b-1}\right)\vspace{2mm}\\
&\le& h\left( \dfrac{2\sqrt{2}}{b-1}\right)+h(d_1b^{l_1}-(d_1-d_2))\vspace{2mm}\\
&\le& h(2\sqrt{2})+h(b-1)+h(d_1)+h(d_1-d_2)+\log 2+l_1h(b)\vspace{2mm}\\
&\le& \dfrac{1}{2}\log 8+3\log 9+\log 2+l_1\log b\vspace{2mm}\\
&\le& \dfrac{1}{2}\log 8+3\log 9+\log 2+\log 9.1+3.86\cdot 10^{13}(1+\log 1.3n)\vspace{2mm}\\
&\le& 3.87\cdot10^{13}(1+\log 1.3n).
\end{array}
\]
Also we have $h(\gamma_2)=\dfrac{1}{2}\log\alpha,\; h(\gamma_3)=\log 10.$ Thus, we can take
\[
A_1:=7.74\cdot 10^{13}(1+\log 1.3n),\quad A_2:=0.89\quad \text{and}\quad A_3:=4.7.
\]
Therefore, we get
\begin{align}\label{log-Gamma_2}
\log|\Gamma_2|>-C_2(1+\log 1.3n)^2,
\end{align}
where $C_2=3.14\times10^{26}.$ By combining \eqref{Gamma_2-bound} and \eqref{log-Gamma_2}, we can get
\[
n\log\alpha<C_2(1+\log 1.3n)^2+\log 3,
\]
this implies that $n<1.82\times 10^{30}.$ Hence we can conclude from \eqref{nlog-alpha} that
$$
l_1+l_2<\dfrac{n\log\alpha+1.5}{\log b}<2.4\times 10^{30}\quad \text{if}\quad b=2
$$
and 
$$
l_1+l_2<\dfrac{n\log\alpha+1.5}{\log b}<1.47\times 10^{30}\quad \text{if}\quad 3\le b\le 10.
$$
We summarize what we have proved so far in the following lemma.
\begin{lemma}\label{l_1+l_2}
If $(n, d_1, d_2, l_1, l_2)$ is a solution in non-negative integers of equation \eqref{main-equation_1}, with $d_1,d_2\in \{0,1,\ldots,9\},$
$d_1\neq d_2$ and $d_1>0,$ then $n<1.82\times 10^{30}.$
Moreover we have, 
$$ 
l_1+l_2<2.4\times 10^{30}\quad \text{if}\quad b=2
$$
and 
$$ 
l_1+l_2<1.47\times 10^{30}\quad \text{if}\quad 3\le b\le 10.
$$
\end{lemma}
\subsection{Reducing the Bound on $n$}

We use the Lemma~\ref{Dujella-Petho} to reduce the bound for $n$. Let
\begin{align*}
\Lambda_1&:=-\log(\Gamma_1+1)\\
          &=(l_1+l_2)\log b-n\log\alpha-\log\left(\dfrac{b-1}{2d_1\sqrt{2}}\right).
\end{align*}
From \eqref{Gamma_1-bound}, we conclude that
\[
\left|e^{-\Lambda_1}-1\right|<\dfrac{9.1}{b^{l_1}}.
\]
Assume that $l_1\ge 5.$ Since $2\le b\le 10,$ we get $\left|e^{-\Lambda_1}-1\right|<\dfrac{9.1}{b^{l_1}}<\dfrac{1}{2},$ which implies that $\dfrac{1}{2}<e^{-\Lambda_1}<\dfrac{3}{2}.$ If $\Lambda_1>0,$ then
$$
0<\Lambda_1<e^{\Lambda_1}-1=e^{\Lambda_1}(1-e^{-\Lambda_1})<\dfrac{18.2}{b^{l_1}}.
$$
If $\Lambda_1<0$, then 
$$
0<|\Lambda_1|<e^{|\Lambda_1|}-1=e^{-\Lambda_1}-1<\dfrac{9.1}{b^{l_1}}.
$$
In any case, it is always holds true $0<|\Lambda_1|<\dfrac{18.2}{b^{l_1}},$ which implies
\begin{align}\label{reduction_1}
0<\left|(l_1+l_2)\dfrac{\log b}{\log\alpha}-n-\dfrac{\log\left((b-1)/2d_1\sqrt{2}\right)}{\log\alpha}\right|<20.7\cdot b^{-l_1}.
\end{align}
Note also that $\dfrac{\log b}{\log\alpha}$ is irrational. In fact, if $\dfrac{\log b}{\log\alpha}=\dfrac{p}{q}$ ($p, q\in \ZZ$ and $p>0, q>0,$ $\gcd(p, q)=1$), then $\alpha^p=b^q\in \ZZ$ which is an absurdity since $2\le b\le 10.$ Taking into account the inequality \eqref{reduction_1} and the Lemma~\ref{Dujella-Petho}, we can take
\[
\tau:=\dfrac{\log b}{\log\alpha},\quad \mu:=-\dfrac{\log\left((b-1)/2d_1\sqrt{2}\right)}{\log\alpha},\quad A:=20.7,\quad B:=b.
\]
According to Lemma~\ref{l_1+l_2}, we can take $M:=2.4\times10^{30}$ for $b=2$ and $M:=1.47\times10^{30}$ for $3\le b\le 10.$ Let $q_t$ be the denominator of the $t$-th convergent of the continued fraction of $\tau.$ We therefore have everything ready to apply Lemma~\ref{Dujella-Petho}. The following table provides information on the results obtained from the applications of Lemma~\ref{Dujella-Petho}.

\begin{tabular}{|c|c|c|c|c|c|c|c|c|c|}
\hline $b$ &2 & 3 & 4 & 5 & 6 & 7 & 8 & 9& 10  \\  \hline 
$q_t$  & $q_{61}$ & $q_{61}$ & $q_{69}$ & $q_{58}$ & $q_{47}$ & $q_{61}$ & $q_{62}$ & $q_{54}$ & $q_{70}$  \\ \hline 
$\varepsilon \ge$ &0.493 & 0.418 & 0.19 & 0.12 & 0.013 & 0.277 & 0.005 & 0.017  & 0.01  \\ \hline
$l_1 \le$ & 112  & 69 & 56 & 47 & 44 & 41 & 40 & 35 & 35 \\ \hline
\end{tabular}\\

It should be noted with regard to the data in the table above that in all cases $1\le l_1\le 112.$ Let
\begin{align*}
\Lambda_2&:=\log(\Gamma_2+1)\\
         &=l_2\log b-n\log\alpha+\log\left(\dfrac{2(d_1b^{l_1}-(d_1-d_2))\sqrt{2}}{b-1} \right).
\end{align*}
From \eqref{Gamma_2-bound} and $n> 110$, we conclude that
\[
\left|e^{\Lambda_2}-1\right|<\dfrac{3}{\alpha^n}<\dfrac{1}{2},
\]
which implies that
$\dfrac{1}{2}<e^{\Lambda_2}<\dfrac{3}{2}.$ If $\Lambda_2>0,$ then $0<\Lambda_2<e^{\Lambda_2}-1<\dfrac{3}{\alpha^n}.$ If $\Lambda_2<0$, then
$$
 0<|\Lambda_2|<e^{|\Lambda_2|}-1=e^{-\Lambda_2}-1=e^{-\Lambda_2}(1-e^{\Lambda_2})<\dfrac{6}{\alpha^n}.
$$ 
It follows in all cases that $0<|\Lambda_2|<\dfrac{6}{\alpha^n},$ thus we have
\begin{align}\label{reduction_2}
0<\left| l_2\dfrac{\log b}{\log\alpha}-n+\dfrac{\log\left(2\sqrt{2}(d_1b^{l_1}-(d_1-d_2))/(b-1) \right)}{\log\alpha}\right|<\dfrac{6.81}{\alpha^n}.
\end{align}
Note that $\dfrac{\log b}{\log \alpha}$ is an irrational number. By referring to \eqref{reduction_2}, we can choose the following data in order to apply Lemma~\ref{Dujella-Petho}.
\[
\tau:=\dfrac{\log b}{\log\alpha},\quad \mu:=\dfrac{\log\left(2\sqrt{2}(d_1b^{l_1}-(d_1-d_2))/(b-1) \right)}{\log\alpha},\quad A:=6.81,\quad B:=\alpha\quad
\]
and 
$M:=2.4\times10^{30}$ for $b=2$ and $M:=1.47\times 10^{30}$ for $3\le b\le 10.$ Let $q_t$ be the denominator of the $t$-th convergent of the continued fraction of $\tau.$ With the help of Maple, we find the following results.

\begin{tabular}{|c|c|c|c|c|c|c|c|c|c|}
\hline $g$ &2 & 3 & 4 & 5 & 6 & 7 & 8 & 9 & 10  \\  \hline 
$q_t$ & $q_{70}$ & $q_{63}$ & $q_{70}$ & $q_{63}$ & $q_{52}$ & $q_{63}$ & $q_{62}$ & $q_{57}$ & $q_{57}$  \\ \hline 
$\varepsilon >$ &0.007 & 0.0009 & 0.001 & 0.0007 & $10^{-5}$ & $2\times10^{-6}$ & 0.0006 & $10^{-4}$& $10^{-4}$   \\ \hline
$n \le$ & 99  & 96 & 94 & 96 & 103 & 103 & 96 & 98& 98 \\ \hline
\end{tabular}\\

Thus, $n\le 103$ is valid in all cases, contradicting the fact that $n>110$. Now, we search for the solutions to the Diophantine equation  \eqref{main-equation_1} with $$
0\le n\le 110,\; 1\le l_1\le 112,\; 1\le l_2\le 142,\; 2\le b\le 10,\; 1\le d_1\le b-1\
$$ 
and $0\le d_2\le b-1,$ by applying a program written in Maple and we only get the solutions listed in Theorem~\ref{main:1}. This completes the proof.

\section{Proof of Theorem~\ref{main:2}}\label{sec:4}
The proof of Theorem~\ref{main:2} is almost similar to that of Theorem~\ref{main:1}, but in this case Legendre's criterion (Lemma~\ref{Legendre}) will be applied in special cases where Lemma~\ref{Dujella-Petho} cannot be applied. To avoid repetitions, we will remove some details in this section.
\subsection{Bounding $n$ and $l_1+l_2$}
According to \eqref{main-equation_2}, we get
\begin{align}\label{Q_n=rep}
Q_n=\dfrac{1}{b-1}\left( d_1 b^{l_1+l_2}-(d_1-d_2)b^{l_2}-d_2\right),
\end{align}
where $d_1, d_2, l_1$ and $l_2$ are non negative integers with $d_1, d_2\in \{1,\ldots,b-1\},$ $d_1\ne d_2$ and $2\le b\le 10.$ 
Throughout this subsection we assume that $n>300.$  From \eqref{P_n-Q_n} and \eqref{Q_n=rep}, we can get 
$$
\alpha^{n-2}\le Q_n<b^{l_1+l_2}\quad \text{and}\quad b^{l_1+l_2-1}<Q_n< \alpha^{n+1},
$$
which implies that
\begin{align*}
(l_1+l_2)\log b-\log\alpha -\log b<n\log\alpha<(l_1+l_2)\log b+2\log \alpha.
\end{align*}
Since $-\log\alpha-\log 10<-\log\alpha-\log b,$ we get what follows 
\begin{align}\label{nlog_alpha}
(l_1+l_2)\log b-3.2<n\log\alpha<(l_1+l_2)\log b+1.8.
\end{align}
 Combining now \eqref{Q_n=rep} and Binet's formula for Pell-Lucas sequence, it is easy to see that
\begin{align*}
\left|(b-1)\alpha^n-d_1b^{l_1+l_2} \right|&\le (b-1)|\beta|^n+|d_1-d_2|b^{l_2}+d_2\\
&\le 3(b-1)b^{l_2}\le 27\cdot b^{l_2},
\end{align*}
which implies that
\begin{align}\label{Gamma_3-bound}
\left| \dfrac{b-1}{d_1}\cdot\alpha^n\cdot b^{-(l_1+l_2)}-1\right|\le \dfrac{27\cdot b^{l_2}}{d_1b^{l_1+l_2}}\le \dfrac{27}{b^{l_1}}.
\end{align}
Let 
\begin{align}\label{Gamma_3}
\Gamma_3:=\dfrac{b-1}{d_1}\cdot\alpha^n\cdot b^{-(l_1+l_2)}-1.
\end{align}
In fact $\Gamma_3\ne 0.$ Indeed if $\Gamma_3=0$, then we would get that
\[
\alpha^n=\dfrac{d_1b^{l_1+l_2}}{b-1}\in \QQ,
\]
which is impossible because $\alpha^n$ is an irrational number for $n\ge 1$.  Thus, we can apply Lemma~\ref{Matveev} on \eqref{Gamma_3} with the data: $t:=3$ and
\[
(\gamma_1,b_1):=\left(\dfrac{b-1}{d_1},1\right),\quad(\gamma_2,b_2):=(\alpha,n),\quad (\gamma_3,b_3):=(b,-l_1-l_2).
\]
Note that $\gamma_1,\; \gamma_2$ and $\gamma_3$ are positive real numbers and elements of the field $\KK=\QQ(\gamma_1,\gamma_2,\gamma_3)=\QQ(\alpha).$ It follows that $D=[\KK:\QQ]=2.$ Moreover, 
\[
h(\gamma_1)=h\left(\dfrac{b-1}{d_1} \right)= \log\left( \max\{b-1,d_1\}\right)\le \log 9,
\]
and
\[ h(\gamma_2)=\dfrac{1}{2}\log\alpha,\quad h(\gamma_3)=\log b\le \log 10.
\]
Thus, according to Lemma~\ref{Matveev} we can take 
$$
A_1=4.4,\quad A_2=0.89\quad\text{and}\quad A_3:=4.7.
$$
From $b^{l_1+l_2-1}<Q_n< \alpha^{n+1}$ and $2\le b\le 10$ with $n> 300$, we easily get that
\[
l_1+l_2<n\dfrac{\log\alpha}{\log 2}+\dfrac{\log2\alpha}{\log 2}<1.3n.
\]
Since $B\ge \max\{1, n, l_1+l_2\},$ we can take $B=1.3n.$ Therefore, we get 
\begin{align}\label{log-Gamma_3}
\log|\Gamma_3|>-C_3(1+\log 1.3n),
\end{align}
where $C_3=1.784\times 10^{13}.$
Hence, from \eqref{Gamma_3-bound} and \eqref{log-Gamma_3}, we have
\begin{align}\label{l_1log b'}
l_1\log b<C_3(1+\log 1.3n)+\log 27.
\end{align}
We transform the equation \eqref{Q_n=rep} again to get something like this
\[
\left|\alpha^n-\left(\dfrac{d_1b^{l_1}-(d_1-d_2)}{b-1} \right)b^{l_2} \right|=\left|\beta^n+\dfrac{d_2}{b-1}\right|\le \alpha^{-n}+\dfrac{d_2}{b-1} <2,
\]
which implies
\begin{align}\label{Gamma_4-bound}
\left|\left( \dfrac{d_1b^{l_1}-(d_1-d_2)}{b-1}\right)\cdot\alpha^{-n}\cdot b^{l_2} -1\right|< \dfrac{2}{\alpha^n}.
\end{align}
Put
\begin{equation}\label{Gamma_4}
\Gamma_4:=\left(\dfrac{d_1g^{l_1}-(d_1-d_2)}{b-1}\right)\cdot\alpha^{-n}\cdot b^{l_2} -1.
\end{equation}
Since assuming $\Gamma_4=0$  leads to
\[
\alpha^n=\left(\dfrac{d_1b^{l_1}-(d_1-d_2)}{b-1}\right)\cdot b^{l_2}\in \QQ,
\]
which is an impossibility, then we must have $\Gamma_4\neq 0.$ 
Thus, we can apply Lemma~\ref{Matveev} on \eqref{Gamma_4} by considering the following data:
\[
t:=3,\quad \gamma_1:=\dfrac{d_1b^{l_1}-(d_1-d_2)}{b-1},\quad \gamma_2:=\alpha,\quad \gamma_3:=b
$$
and
$$
b_1:=1,\quad b_2:=-n,\quad b_3:=l_2.
\]
Note also that $\gamma_1,\; \gamma_2$ and $\gamma_3$ are positive real numbers and elements of the field $\KK=\QQ(\gamma_1,\gamma_2,\gamma_3)=\QQ(\alpha).$ So, we have $ D=[\KK:\QQ]=2.$ Furthermore, from \eqref{l_1log b'} we get
\begin{align*}
h(\gamma_1)&=h\left(\dfrac{d_1b^{l_1}-(d_1-d_2)}{b-1}\right)\\
&\le h\left( d_1b^{l_1}-(d_1-d_2)\right)+h(b-1)\\
&\le 3\log (b-1)+l_1\log b+\log2\\
&\le 2\times10^{13}(1+\log 1.3n).
\end{align*}
Thus, we can take 
$$
A_1=4\times10^{13}(1+\log 1.3n)),\quad  A_2=0.89\quad\text{and}\quad A_3=4.7.
$$
Using $l_1+l_2<1.3n$, we can take $B:=1.3n.$ Hence, Lemma~\ref{Matveev} tells us that 
\begin{align}\label{log-Gamma_4}
\log|\Gamma_4|>-C_4(1+\log 1.3n)^2,
\end{align}
where $C_4=1.63\times10^{26}.$ By combining \eqref{Gamma_4-bound} and \eqref{log-Gamma_4}, we can get
\[
n\log\alpha<C_4(1+\log 1.3n)^2+\log 2,
\]
this implies that $n<9.2\times 10^{29}.$ It follows from \eqref{nlog_alpha} that
$$
l_1+l_2<\dfrac{n\log\alpha+3.2}{\log b}<1.17\times 10^{30}\quad \text{if}\quad b=2 
$$
and 
$$
l_1+l_2<\dfrac{n\log\alpha+3.2}{\log b}<7.39\times 10^{29}\quad \text{if}\quad 3\le b\le 10.
$$
In summary we have the following result.
\begin{lemma}\label{l_1+l_2'}
If $(n, d_1, d_2, l_1, l_2)$ is a solution in non-negative integers of equation \eqref{main-equation_2}, with $d_1,d_2\in \{0,1,\ldots,9\},$
$d_1\neq d_2$ and $d_1>0,$ then $n<9.2\times 10^{29}.$
Moreover we have, 
$$
l_1+l_2<\dfrac{n\log\alpha+3.2}{\log b}<1.17\times 10^{30}\quad \text{if}\quad b=2 
$$
and 
$$
l_1+l_2<\dfrac{n\log\alpha+3.2}{\log b}<7.39\times 10^{29}\quad \text{if}\quad 3\le b\le 10.
$$
\end{lemma}

\subsection{Reducing the Bound on $n$}

We use Lemmas~\ref{Dujella-Petho} and \ref{Legendre} to reduce the bound for $n.$ Put
 \begin{align*}
 \Lambda_3&:=-\log(\Gamma_3+1)\\
          &=(l_1+l_2)\log b-n\log\alpha-\log\left(\dfrac{b-1}{d_1}\right).
 \end{align*}
From \eqref{Gamma_3-bound}, we have
\begin{align}
\left| e^{-\Lambda_3}-1\right|<\dfrac{27}{b^{l_1}}.
\end{align}
If $l_1\ge 6,$ then $|e^{-\Lambda_3}-1|<\dfrac{27}{b^{l_1}}<\dfrac{1}{2},$ which implies that 
 $$
 0<|\Lambda_3|<\dfrac{54}{b^{l_1}}.
 $$
 Thus
\begin{align}\label{red}
0<\left|(l_1+l_2)\dfrac{\log b}{\log\alpha}-n-\dfrac{\log((b-1)/d_1)}{\log\alpha}\right|<62\cdot b^{-l_1}.
\end{align}
In fact, we need to see the following two cases.

\textbf{Case $d_1\ne b-1.$}

 According to \eqref{red} and Lemmas~\ref{Dujella-Petho} and \ref{l_1+l_2'}, we can  take $M:=7.39\cdot 10^{29}$ if $3\le b\le 10.$ Towards applying Lemma~\ref{Dujella-Petho} for $3\le b\le 10$ and $1\le d_1\le b-2$, we define the following quantities
\[
\tau:=\dfrac{\log b}{\log\alpha},\quad \mu:=-\dfrac{\log((b-1)/d_1)}{\log\alpha},\quad A:=62,\quad B:=b.
\]
Also, it is easy to see that $\dfrac{\log b}{\log\alpha}$ is an irrational number. Let $q_t$ be the denominator of the $t$-th convergent of the continued fraction of $\tau.$ The results obtained following the application of the Lemma~\ref{Dujella-Petho} are presented as can be seen in the following table
$$
\begin{tabular}{|c|c|c|c|c|c|c|c|c|c|}
\hline $b$  & 3 & 4 & 5 & 6 & 7 & 8 & 9 & 10 \\  \hline 
$q_t$  & $q_{59}$ & $q_{67}$ & $q_{58}$ & $q_{45}$ & $q_{60}$ & $q_{61}$ & $q_{54}$ & $q_{68}$  \\ \hline 
$\varepsilon >  $  & 0.26 & 0.16 & 0.19 & 0.01 & 0.09 & 0.09 & 0.08& 0.006   \\ \hline
$l_1 \le $  & 69 & 55 & 47 & 44 & 39 & 38 & 35 & 35 \\ \hline
\end{tabular}
$$
It follows that,
\begin{align}\label{l_1}
l_1\le \dfrac{\log \left(62q_t/\varepsilon\right)}{\log b}\le 69,
\end{align}
which holds in all cases.  

\textbf{Case $d_1=b-1.$}

In this case we need to apply Lemma~\ref{Legendre} since $\mu =0.$ The inequality \eqref{red} can be rewritten as
\begin{align*}
0<\left|(l_1+l_2)\dfrac{\log b}{\log\alpha}-n\right|<\dfrac{62}{b^{l_1}}.
\end{align*}
Referring to Lemmas~\ref{Legendre} and \ref{l_1+l_2'}, we can take $M:=1.17\times 10^{30}$ if $b=2$ and $M:=7.39\times10^{29}$ if $3\le b\le 10.$ For $2\le b\le 10,$ we use Maple to find the first convergent $q_N$ such that $q_N>M$ and then we get $a(M):=\max\{a_i:i=0,\ldots,N\}.$ Therefore, Lemma~\ref{Legendre} tells us that 
\begin{align*}
\dfrac{62}{b^{l_1}}>\left|(l_1+l_2)\dfrac{\log b}{\log\alpha}-n\right|>\dfrac{1}{(a(M)+2)(l_1+l_2)},
\end{align*}
which implies 
\[
l_1<\dfrac{\log\left(62\cdot(a(M)+2)\cdot (l_1+l_2) \right)}{\log b}.
\]
Thus we obtain the following results which follow from the application of Lemma~\ref{Legendre}.
$$
\begin{tabular}{|c|c|c|c|c|c|c|c|c|c|}
\hline $b$ &2 & 3 & 4 & 5 & 6 & 7 & 8 & 9 & 10 \\  \hline 
$q_N>M$ &$q_{59}$  & $q_{58}$ & $q_{67}$ & $q_{56}$ & $q_{44}$ & $q_{59}$ & $q_{58}$ & $q_{52}$& $q_{67}$   \\ \hline 
$a(M)$ &100 & 130 & 110 & 163 & 509 & 33 & 34 & 68 & 52  \\ \hline
$l_1 \le $ &112 & 70 & 55 & 48 & 44 & 39 & 36 & 35 & 34\\ \hline
\end{tabular}
$$
So we have 
\begin{align}\label{l_1'}
l_1\le 112.
\end{align}
By combining \eqref{l_1} and \eqref{l_1'}, we see that $1\le l_1\le 112$ holds in all cases. 

Put now 
\begin{align*}
\Lambda_4&:=\log(\Gamma_4+1)\\
         &=l_2\log b-n\log\alpha+\log\left(\dfrac{d_1b^{l_1}-(d_1-d_2)}{b-1} \right).
\end{align*}
 Since $n> 300$, we can conclude from \eqref{Gamma_4-bound} that 
\[
\left|e^{\Lambda_4}-1\right|<\dfrac{2}{\alpha^n}<\dfrac{1}{2},
\]
which implies that $0<|\Lambda_4|<\dfrac{4}{\alpha^n}$ and therefore
\begin{align}\label{red-4}
0<\left| l_2\dfrac{\log b}{\log\alpha}-n+\dfrac{\log\left((d_1b^{l_1}-(d_1-d_2))/(b-1) \right)}{\log\alpha}\right|<4.6\cdot \alpha^{-n}.
\end{align}
It is necessary to specify that the case $b=2$ is only possible if $d_1=1$ and $d_2\in\{0,1\}$. So we need to study it in a special way. 

$\bullet$ If $d_2=0$ and $l_1\neq 1$, then \eqref{red-4} becomes
\begin{equation}
0<\left| l_2\dfrac{\log 2}{\log\alpha}-n+\dfrac{\log\left(2^{l_1}-1 \right)}{\log\alpha}\right|<4.6\cdot \alpha^{-n}.
\end{equation}
So, in this case we apply Lemma~\ref{Dujella-Petho} with the data:
\[
\tau=\dfrac{\log 2}{\log\alpha},\quad \mu:=\dfrac{\log\left(2^{l_1}-1 \right)}{\log\alpha},\quad A:=4.6,\quad B:=\alpha.
\]
Also, we can take $M:=1.17\cdot 10^{30}$. Using Maple, we find that the denominator $q_{62}$ of the $62$-th convergent of the continued fraction of $\log 2/\log \alpha$ satisfies $q_{62}>6M$ and $\varepsilon> 0.00175$. So it follows from Lemma~\ref{Dujella-Petho} that 
\begin{equation}\label{1}
n<\dfrac{\log(4.6 q_{62}/0.00175)}{\log \alpha}<93.
\end{equation}
$\bullet$ If $d_2=0$ and $l_1=1$ or $d_2=1,$ then the relation \eqref{red-4} becomes
\begin{equation}\label{a1}
0<\left| \lambda\dfrac{\log 2}{\log\alpha}-n\right|<4.6\cdot \alpha^{-n}, \quad \text{where } \lambda\in\{l_2,l_1+l_2\}.
\end{equation}
Since the denominator  $q_{59}$ of the $59$-th convergent of the continued fraction of $\log 2/\log \alpha$ satisfies $q_{59}>M$ and $a(M)=100$, then by Lemma~\ref{Legendre}, we get
\begin{equation}\label{a2}
\left| \lambda\dfrac{\log 2}{\log\alpha}-n\right|> \dfrac{1}{102\lambda}>\dfrac{1}{102\cdot 1.17\cdot 10^{30}}.
\end{equation}
Combining \eqref{a1} and \eqref{a2}, we deduce that
\begin{equation}\label{2}
n<\dfrac{\log(4.6 \cdot 102\cdot 1.17 \cdot 10^{30})}{\log \alpha}< 85.
\end{equation}
From now on we will see what happens with $3\le b \le 10$. 
For this, we study the following two cases while exploiting the inequality \eqref{red-4}.

\textbf{Case $(d_1,l_1,d_2)\neq (1,1,0).$} 

Note that in this case, by referring to the inequality \eqref{red-4} we are able to apply Lemma~\ref{Dujella-Petho} while choosing the following data
\[
\tau:=\dfrac{\log b}{\log\alpha},\quad \mu:=\dfrac{\log\left((d_1b^{l_1}-(d_1-d_2))/(b-1) \right)}{\log\alpha},\quad A:=4.6,\quad B:=\alpha
\]
and 
$M:=7.39\times10^{29}.$ Let $q_t$ be the denominator of the $t$-th convergent of the continued fraction of $\tau.$ With the help of Maple, we get the following results.
$$
\begin{tabular}{|c|c|c|c|c|c|c|c|c|}
\hline $b$ & 3 & 4 & 5 & 6 & 7 & 8 & 9 & 10 \\  \hline 
$q_t$ & $q_{81}$ & $q_{102}$ & $q_{102}$ & $q_{100}$ & $q_{123}$ & $q_{120}$ & $q_{134}$  &$q_{145}$ \\ \hline 
$\varepsilon >$ & $10^{-12}$ & $10^{-17}$ &$10^{-26}$ & $10^{-30}$ & $10^{-33}$ & $10^{-29}$ & $10^{-38}$& $10^{-43}$   \\ \hline
$n\le$ & 142 & 178 & 208 & 232 & 250 & 268 & 283& 297 \\ \hline
\end{tabular}
$$
So, we have in all cases 
\begin{align}\label{bound-n-1}
n\le 297.
\end{align}
\textbf{Case $(d_1,l_1,d_2)=(1,1,0).$} 

Here the relation \eqref{red-4} becomes
\begin{align}\label{red-4-Legender-last}
0<\left| l_2\dfrac{\log b}{\log\alpha}-n\right|<\dfrac{4.6}{\alpha^n}.
\end{align}
Next, we apply Lemma~\ref{Legendre} with $M:=7.39\cdot 10^{29}$ while finding $q_N$ such that $q_N>M$ and $a(M):=\{a_i:i=0,1,\ldots,N\}$. We have 
\begin{align}\label{red-4-Legender'-last}
\left|l_2\dfrac{\log b}{\log\alpha}-n\right|>\dfrac{1}{(a(M)+2)\cdot l_2}>\dfrac{1}{(a(M)+2)\cdot 7.39\times 10^{29}}.
\end{align}
By referring to \eqref{red-4-Legender-last} and \eqref{red-4-Legender'-last}, we obtain 
\[
n<\dfrac{\log\left(4.6\cdot (a(M)+2)\cdot 7.39\cdot 10^{29}\right)}{\log\alpha}.
\]
According to the above equality and $3\le b\le 10$ , we get the results as follows  thanks to Maple.
$$
\begin{tabular}{|c|c|c|c|c|c|c|c|c|c|}
\hline $b$  & 3 & 4 & 5 & 6 & 7 & 8 & 9 & 10  \\  \hline 
$q_N>M$  & $q_{58}$ & $q_{67}$ & $q_{56}$ & $q_{44}$ & $q_{59}$ & $q_{58}$ & $q_{52}$ & $q_{67}$  \\ \hline 
$a(M)$  & 130 & 110 & 163 & 509 & 33 & 34 & 68& 52   \\ \hline
$n\le$  & 85 & 85 & 85 & 86 & 83 & 83 & 84& 84 \\ \hline
\end{tabular}
$$
It follows that
\begin{align}\label{bound-n-2}
n\le 86.
\end{align}
From the relations \eqref{1}, \eqref{2}, \eqref{bound-n-1} and \eqref{bound-n-2}, we easily conclude that $n\le 297$. This contradicts the assumption $n>300.$ Finally, we search for the solutions to the Diophantine equation \eqref{main-equation_2} with 
$$
0\le n\le 300,\; 1\le l_1\le 112,\; 1\le l_2\le 386,
$$
and
$$
2\le b\le 10,\; 1\le d_1\le b-1,\; 0\le d_2\le b-1,
$$ 
by applying a program written in Maple and we only get the solutions listed in Theorem~\ref{main:2}. This completes the proof.

\section*{Acknowledgements} The first author is supported by Institut de Math\'ematiques et de Sciences Physiques de l'Universit\'e d’Abomey-Calavi. The second author is supported by the Croatian Science Fund, grant HRZZ-IP-2018-01-1313.

%
Institut de Math\'ematiques et de Sciences Physiques, Universit\'e d'Abomey-Calavi, B\'enin \\
Email: adedjnorb1988@gmail.com \\[6pt]
Faculty of Civil Engineering, University of Zagreb, Fra Andrije
Ka\v{c}i\'{c}a-Mio\v{s}i\'{c}a 26, 10000 Zagreb, Croatia \\
Email: filipin@grad.hr \\[6pt]
Department of Mathematics, Institute of Science and Technology, University Center of Mila, Algeria\\
Email salahrihane@hotmail.fr\\[6pt]
Department of Mathematics, Statistics and Computer Science, Purdue
University Northwest, 1401 S, U.S. 421, Westville IN 46391 USA \\
Email: atogbe@pnw.edu\\[6pt]

\end{document}